\documentclass[letterpaper, 10pt, conference]{ieeeconf} 
\IEEEoverridecommandlockouts 
\usepackage{latin1}
\usepackage{amsfonts}
\IEEEoverridecommandlockouts                              
\usepackage{graphics} 
\usepackage{epsfig} 

\def\ip#1#2{\left\langle #1, #2 \right\rangle }
\def\triangleeq{\stackrel{\triangle }{=}}
\def\E{\rm E}
\def\trace{{\rm tr}}
\def\psd{{\sc psd}}
\def\dis{d}
\def\DIS{D}

\def\cP{\mathcal P}
\def\cQ{\mathcal Q}
\def\cR{\mathcal R}
\def\cD{\mathcal D}

\def\b{\mathbf b}
\def\c{\mathbf c}

\def\q{\mathbf q}
\def\r{\mathbf r}

\def\Fr{{\mathfrak F}_\r}
\def\Fs{{\mathfrak F}_\Sigma}

\newtheorem{example}{Example}[section]

\title{\LARGE \bf {Approximative Covariance interpolation}}

\author{Per Enqvist
\thanks{This work was supported by the Swedish research council.}
\thanks{The author is with the division of Optimization and Systems Theory, Department of Mathematics,
        Royal Institute of Technology, SE-100 44 Stockholm, Sweden.
        {\tt\small penqvist@math.kth.se}}%
}

\begin{document}

\maketitle
\thispagestyle{empty}
\pagestyle{empty}

\begin{abstract}
When methods of moments are used for identification of power spectral
densities, a model is matched to estimated second order statistics
such as, {\it e.g.}, covariance estimates.  If the estimates are good
there is an infinite
family of power spectra consistent with such an estimate and in
applications, such as identification, we want to single out the most
representative spectrum. 
We choose a prior spectral density to
represent {\it a priori} information, and the spectrum closest to it
in a given quasi-distance is determined.
However, if the estimates are based on few data, or the model class 
considered is not consistent with the process considered, it may 
be necessary to use an approximative covariance interpolation. 
Two different types of regularizations are considered in this paper
that can be applied on many covariance interpolation based 
estimation methods.

\end{abstract}

\section{Introduction}

Most system identification methods are based on an 
algorithm that is proven to give efficients estimates when 
the number of data goes to infinity. 
One such common estimate is the maximum likelihood method. 
However, in many cases only a small amount of data is available and 
the estimation method may give unexpected results. 
Here we will consider methods based on covariance interpolation instead. 
Depending on which model class is considered there are a number of different 
methods around now for matching AR, MA, ARMA and other models to 
covariances, such as the ones derived by Lindquist, Byrnes, Georgiou, Pavon, Ferrante, 
{\it et. al.} based on minimizing the
Kullback-Leibler \cite{GL}, Hellinger \cite{HellingerMIMO},  
the Itakura-Saito quasi-distance \cite{Shore,ShoreJohnson,EK}, 
and other distance concepts.
However, also these methods depends on the amount of data that is available
and also structural constraints. 
The covariances have to be estimated from the data and the errors in the 
estimates will increase the smaller the available data set is. 
Estimating the covariances from a short data sequence may generate a 
covariance matrix that is not non-negative definite, or does not have a 
supposed Toeplitz structure or the estimate does not 
correspond to a spectra in the supposed model class.
So for short data sequences it is necessary to regularize the methods to obtain
relevant model estimates. 
In this paper we compare two different approaches 
for dealing with these kinds of problems; the two different kinds of 
regularizations are  based on
 quadratic penalties on the covariance estimation errors  and 
extra entropy regularization of the determined spectrum.
These approaches have been used before for the maximum entropy 
method for AR-models, the Kullback-Leibler method for 
 the ARMA case with fixed MA-part and a combined covariance
and cepstrum interpolation problem, but here they will 
be used and compared in a more general setting.

The first kind of problem, with non-negative definite covariance matrices, 
is often ``solved'' by using a biased estimate of the covariance matrix.
This bias is usually small and goes to zero as the number of data 
grows, but for small data sets it can be relevant. 
Another approach is to use a regularization of the first kind mentioned above, i.e., 
to find a spectrum within the model class which has a 
small quadratic distance to the estimated covariance matrix. 
By combining the covariance interpolation methods based on entropy maximization
with a quadratic distance penalty the structure of the spectrum is taken into 
account when the best covariance sequence close to the estimates is determined.

The second kind of problem, with the estimate of the covariance matrix not having 
the supposed structure, is often solved using a projection onto the class of 
matrices with the desired structure. 
This problem is most obvious when a state-covariance interpolation approach is used;
Then there is an imposed structure determined by the $(A,b)$ matrices in the 
state-covariance definition.
Again, another approach is to use the regularization of the first kind mentioned 
above. A small distance to a matrix with the desired structure is then obtained.

The third kind of problem, with a covariance estimate that can not be interpolated
by a spectrum in the model class (but has the desired structure and is non-negative) 
as in MA-model covariance interpolation for some covariance estimates.
Probably the most common approach to resolve this problem 
 is to project the covariance estimates 
onto the set of covariances feasible for the desired model class.
For the MA case, this would be the projection onto a positive cone, but to avoid 
having zeros on the unit circle a projection to a slightly smaller cone 
should be performed. 
Another approach is to use a regularization of either the first or second 
kind mentioned above. The amount of quadratic penalty regularization for the 
first method has to be determined recursively, and might fail for some 
cases as will be shown by some examples. 
The extra entropy regularization treats this case in an easier way and finds both 
the best approximating valid covariance and the interpolant with one 
optimization problem. 

If we want to determine a MA-model estimate for a state-covariance
estimated from a short data sequence, 
all of the three kinds of problems described above may occur.
Then it would be necessary to use a combination of the 
two types of regularizations.



Here, a general approach to the covariance matching problem is 
taken that holds for a large set of different quasi-distances and 
is inspired by the work in \cite{PerECC09}.


\section{Background}\label{sec:back}

Let $(\ldots,\, y_{-1},\, y_0,\, y_1,\,\ldots)$ be a scalar 
stationary stochastic real valued mean-zero process with covariances
$r_k = \E \{ y_{\ell+k} y_{\ell} \}$ and \psd \ $\Phi$.
The power spectral density $\Phi$ represents the energy content
of the process across frequencies and has the covariances as
Fourier coefficients,
\[
\Phi(e^{i\theta}) \triangleeq \sum_{k=-\infty}^\infty r_k e^{ik\theta}.
\]

Consider the Hilbert space $L_2(-\pi,\pi]$
with the inner product
\[
\ip{a}{b} =
\frac{1}{2\pi}\int_{-\pi}^\pi a(e^{i\theta})b(e^{-i\theta})d\theta.
\]
Then the covariances are given by
$r_k 
= \ip{\Phi}{z^k}. $
Given a finite window of covariances
$\r= \left( \begin{array}{ccccc}
r_0 & r_1 & \ldots & r_n
\end{array}\right)$, 
let $\Fr$ denote the set of \psd \ 
consistent with $\r$, {\it i.e.},
\[
\Fr = \left\{  \Phi \geq 0  \; |  
 \ip{\Phi}{z^k}=r_k, \quad k=0,1,\cdots,n \right\}.
\]
In this paper, $\Phi\geq 0$ means that this inequality should hold on the unit circle,
{\it i.e.}, $\Phi(e^{i\theta}) \geq 0$ for $\theta\in(-\pi,\pi]$.

Furthermore, we assume initially that the 
symmetric  Toeplitz matrix
of the covariances $\r$,
\begin{equation}\label{toep}
 T( \r)  = \left[ \begin{array}{cccc}
 r_0 &  r_1 & \cdots &  r_n \\
 r_1 &  r_0 & \ddots & \vdots \\
\vdots & \ddots & \ddots &  r_1 \\
 r_n & \cdots &  r_1 &   r_0
\end{array} \right]
\end{equation}
is positive definite, hence
the set $\Fr$ contains an infinite number of \psd s \ 
\cite[Sec.6.5]{Porat}.
Let $ \cR = \left\{  \r  \; |  T(\r) > 0 \right\}$.


\section{Moment Matching}

In many situations it is desired to fit a spectral density to data
by finding one of a particular structure by matching moments. 
The most common \psd \ used to model stationary stochastic processes 
are the ones that correspond to Moving-Average (MA) and Auto-Regressive (AR)
processes.
Assume that $Q(z)$ is a pseudo-polynomial of degree $n$, {\it i.e.},
\begin{equation}\label{Q}
Q(z)  = q_0 + \frac{1}{2}q_1(z+z^{-1}) + \cdots +\frac{1}{2}q_n(z^n+z^{-n}).
\end{equation}
Then, $\Phi = Q$ is the \psd \ of a MA-process and $\Phi=1/Q$ is the
\psd \ of an AR-process.  It is well known that for an AR-process the
coefficients $\{ q_k\}_{k=0}^n$ of $Q$ can always be tuned so that a
window of covariances 
$\r \in \cR$ 
is matched.  On the
other hand, it is also well known that for an MA-process there are
some $\r\in\cR$ (actually open subsets of such covariances) 
that are not matched for any choice
of coefficients $\{q_k\}_{k=0}^n$.  In both cases there are $n+1$
parameters that should be tuned to match $n+1$ constraints, but it is
clearly the structure of the \psd \ that determines if solutions
exists or not.

To generalize, let 
\begin{equation}
G(z) := (I-zA)^{-1}B.
\end{equation}
define an input-to-state map, where 
we will assume that $A$ is an 
$n\times n$-stability matrix, $b$ is an 
$n\times 1$ vector and 
 $(A,B)$ is a reachable pair.
Then, for a symmetric matrix  we can define a 
generalized pseudopolynomial 
\begin{equation}
Q(z) := G(z)\Lambda G^*(z).
\end{equation}

Similarly, we generalize $\Fr$ to
\[
\Fs = \left\{  \Phi \geq 0  \; \left|  
 \int G\Phi G^*= \Sigma \right.
 \right\}.
\]

To evaluate the properties of different \psd \ structures, we let 
$\Phi$ depend on $Q$, and it will also be
allowed to depend on some ``prior estimate'' \psd \ $\Psi$.  Assuming
now that $\Phi = F(Q,\Psi)$, the moment matching constraint $\Phi
\in \Fs$, can be expressed as
\begin{equation}\label{momcon}
\int G F(Q,\Psi)G^* =\Sigma =\int G R G^*,
\end{equation}
where $R$ is an arbitrary function in $\Fs$. 


\section{Exact and Approximative interpolation}

If we use unbiased estimates of the state-covariances
from  a realization with a \psd \ $\Phi$ we know that  
the \psd determined by exact moment matching will converge to $\Phi$
as the number of samples tend to infinity, 
if $\Phi$  is in the class of spectrums considered. 

For short realizations it may be necessary to introduce some 
bias to get reasonable estimates. 
By introducing bias the variance of the estimates can be reduced.
How this is done is an important issue.

\subsection{Exact interpolation}

The distance measure will be assumed to be differentiable
in the first argument, and it will be assumed to be a 
quasi-distance, {\it i.e.}, it is assumed that 
$\DIS(\Phi||\Psi)\geq 0$ and $\DIS(\Psi||\Psi)=0$ for any pair of \psd \ $\Phi$ 
and $\Psi$.
Furthermore, we assume that 
\[ \DIS(\Phi||\Psi) = 
\int \dis(\Phi||\Psi). 
\]

Note that $\DIS$ is not assumed to be symmetric, convex, to satisfy the triangle
inequality or be zero if and only if $\Phi=\Psi$.
However, these are certainly desired properties.
Consider the optimization problem, to minimize the 
distance to $\Psi$ for all  $\Phi \in \Fs$, {\it i.e.}:
\begin{equation}\label{exactprimal}
(\cP_=) \quad \left[
\begin{array}{rl}
{\displaystyle \inf_{\Phi\geq 0}} &  {\displaystyle \DIS(\Phi||\Psi) } \\[2ex]
\mbox{s.t.} & 
\int G \Phi G^* - \Sigma =0.
\end{array}
\right]
\end{equation}
Note that here that the \psd \ $\Phi$
 is not constrained to be of a certain form,
this form will be determined by the optimality conditions of the 
Lagrange relaxed functional, which in turn is determined by the
geometry imposed by the distance measure.

The optimization problem $(\cP_=)$ has no finite dimensional
parametrization, but by considering the dual, 
an optimization problem with a finite number of variables is obtained.
To this end, formal calculations are performed to determine the dual.

Form the Lagrangian function
\[
L_0(\Phi;\q) \triangleeq \DIS(\Phi||\Psi)+
\trace \left\{ \Lambda (\Sigma - \int G\Phi G^*)\right\}
\]
and since $\Phi$ is symmetric
 $\ip{\Phi}{\sum_{k=0}^n q_kz_k} = \ip{\Phi}{Q(z)}$, where $Q$ is defined in (\ref{Q}).
Let $R\in \Fs$ arbitrary.
Then the Lagrangian function can be written as
\[
L_0(\Phi; Q) = \DIS(\Phi||\Psi)
+ \trace \{\Lambda\Sigma\} - \ip{\Phi}{Q}.
\]

Assuming that a minimizer exists
let 
\begin{equation}\label{L0minimizer}
\hat \Phi := \arg\min_{\Phi \geq 0} L_0(\Phi,Q),
\end{equation}
this defines the optimal \psd \ as a function of $Q$, {\it i.e.},
\begin{equation}\label{F}
\hat \Phi = F(Q;\Psi) 
\end{equation}
and determines the dual objective function
\begin{equation}
\Omega_0 (Q;\Psi)
 \triangleeq L_0(\hat \Phi,Q)
=L_0(F(Q;\Psi),Q).
\end{equation}

To ensure that the spectral densities $F(Q;\Psi)$ are non-negative
the domain $\cQ$ of feasible $Q$ has to be specified, 
{\it i.e.}, 
\[ \cQ = \{ Q \, | \, F(Q; \Psi)  \geq 0 \}.  \]

This leads to the dual problem to 
determine the maximizer of $\Omega_0$ over all $Q\in\cQ$,
{\it i.e.}, 
\begin{equation}  (\cD_=) \quad 
\left[
\begin{array}{rl}
{\displaystyle \sup_{Q}} & \Omega_0(Q;\Psi,R)\\
\mbox{s.t.} & F(Q;\Psi) \geq 0.
\end{array}
\right],
\end{equation}

The derivative of $\Omega$ is 
(compare the proof of Proposition 4.1 in 
\cite{PerECC09})
\[ \frac{\partial \Omega_0}{\partial Q} =
\int \left( \frac{\partial }{\partial Q} d(F||\Psi)-Q
 \right)\frac{\partial F}{\partial Q}  +
\int \left( R-F
 \right)  
 \]
and using that $F(Q,\Psi)$ minimizes $L_0$
it can be shown that the first integral is zero.
The stationarity conditions for $(\cD_=)$ are then
\[ \Sigma -\int G F(Q;\Psi) G^* =0, \] 
for $k=0,1,\cdots,n$, which ensures that 
for an interior point solution the optimal 
$\Phi\in\Fs$.

When the state-covariance $\Sigma$ is estimated from 
a short sequence of data, it is quite likely that the 
there will exist no exact interpolants. Even for long
data sequences the existence of solutions may fail  
if the given realization does not match  
the class of \psd s considered.

\subsection{Primal regularization}

Consider now the approximative interpolation problem:
\begin{equation}\label{approxprimal}
(\cP_\approx^2) \quad \left[
\begin{array}{rl}
{\displaystyle \inf_{\Phi\geq 0}} &  {\displaystyle
\DIS(\Phi||\Psi) } + \trace\{ DWD \}\\[2ex]
\mbox{s.t.} & 
\int G \Phi G^* - \Sigma = D
\end{array}
\right]
\end{equation}
In this problem we consider not only \psd s  in $\Fs$, but
any \psd \ and then we penalize deviations from the 
nominal state-covariance $\Sigma$ using a quadratic 
penalty term. 
An alternative approach would be to make a fixed extension of the 
set $\Fs$, such as fixed intervals of the parameters in $\Sigma$, that 
approach is taken in \cite{BLunc}.
Once again, the spectral density is not constrained to be of a certain form,
this form will be determined by the optimality conditions of the 
Lagrange relaxed functional, which in turn is determined by the
geometry imposed by the distance measure.


We show that the structure of the optimal $\Phi$ will be the same 
as for $(\cP_=)$.
Form the Lagrangian function
\[
L(\Phi, D;\Lambda) \triangleeq \DIS(\Phi||\Psi) 
+\trace\{ DWD \} \qquad
\]
\[ \quad +\trace\left\{ \Lambda (D+\Sigma - \int G\Phi G^*) 
\right\}
\]
\[
\quad = L_0(\Phi;\q)  
+\trace\{ DWD\} + \trace\{\Lambda D\}
\]


The optimal $\hat D = -\frac{1}{2}W^{-1}\Lambda$.
So for large $W$ the approximation 
errors go to zero (if an exact solution exists).

The optimal \psd \ $\Phi$ is again determined 
by (\ref{L0minimizer}), hence the structure of $\Phi$ 
is preserved and 
 $\hat \Phi = F(Q;\Psi)$, see (\ref{F}).

The dual objective function is then given by 
\begin{equation}
\Omega (Q;\Psi)
=L(\hat \Phi,\hat \Delta,Q)
=\Omega_0 (Q;\Psi)-\frac{1}{4}\trace\{\Lambda W^{-1}\Lambda\}.
\end{equation}

This leads to a dual problem on the form

\begin{equation}  (\cD_\approx^2) \quad 
\left[
\begin{array}{rl}
{\displaystyle \sup_{Q}} & \Omega(Q;\Psi)\\
\mbox{s.t.} & F(Q;\Psi) \geq 0.
\end{array}
\right].
\end{equation}

The stationarity conditions for $(\cD_\approx^2)$ are 
\[  \Sigma - \int G F(Q;\Psi)G^* = \frac{1}{4} \left(
\Lambda W^{-1}+W^{-1} \Lambda \right).
\]

\subsection{Dual regularization}



Consider the dual regularized optimization problem:

\begin{equation}  (\cD_\approx^1) \quad 
\left[
\begin{array}{rl}
{\displaystyle \sup_{Q}} & \Omega_0(Q;\Psi,R)+
\lambda B(Q)\\
\mbox{s.t.} & F(Q;\Psi) \geq 0.
\end{array}
\right],
\end{equation}
where 
$B(Q)$ is a barrier type of function
who´s purpose is to keep the optimum in an 
interior point, and regularize the solution,
{\it i.e.} avoid too sharp pikes in the \psd.

The barrier function will typically be a function like
\[ B_1(Q) =  \int \log (1+Q), \]
whose derivative in the direction of the boundary 
goes to infinity as $\Lambda$ goes to the 
boundary, or 
\[ B_2(Q) = 1- \int \frac{1}{1+Q}, \]
whose function values goes to infinty at the boundary.

The stationarity conditions are then  
\[  \Sigma - \int G F(Q;\Psi)G^* = \lambda \int G
\frac{1}{1+G^*\Lambda G} G^*
\]
and
\[  \Sigma - \int G F(Q;\Psi)G^* = \lambda \int G
\frac{1}{(1+G^*\Lambda G)^2} G^*
\]
respectively.


The right hand side will be small for small $\lambda$. 
If $Q$ is close to zero for some frequencies, the 
integral will still be bounded but have a derivative 
that goes to infinity as $Q$ goes to zero.

The problem $  (\cD_\approx^1) $ is a 
convex optimization problem and could therefore 
be the dual of some optimization problem, but 
the author has not been succesful in finding 
such a primal problem.
For some cases, for example when $\Phi=\Psi/Q$, 
the extra term in the objective function can 
be seen to increase the entropy of the 
resulting \psd.

\subsection{Comparison of the two regularizations}

We note that both the regularizations results in adding 
a concave function of $Q$ to the dual objective function.
In $ (\cD_\approx^1) $ it is a logaritmic term that 
works as a barrier function making sure that the 
optimum is in an interior point of $\cQ$.
If the optimum of the primal problem $(\cP_=)$
is in an interior point, the regularization term 
is rather small and does not affect the solution much
but tends to pull it slightly towards a 
spectrum with \psd \ $g(\Psi)$.
If the optimum of the primal problem $(\cP_=)$
is on the boundary, the unbounded derivative of the 
regularization term will push the solution towards the 
interior.

In  $(\cD_\approx^2)$ the regularization term is a quadratic 
function of the matching error. By allowing a slack in the 
covariance matching constraint the distance $D(\Phi||\Psi)$
can be made smaller and a \psd \ closer to the prior
is obtained. This means that more trust is put on the 
prior information and less is put on the covariances, 
which makes sense if the covariances are estimated from 
short data sequences.
For the Kullback-Leibler distance it is shown in 
\cite{EAtech} that even if the covariances are not in $\cR$,
{\it i.e. }, corresponds to a positive definite Toeplitz matrix, 
an approximative solution is obtained if the $\alpha$  is chosen 
small enough. Note that the covariances can fail to correspond
to a positive definite matrix and they can also fail to form 
a Toeplitz matrix, but an approximation is anyway guarranteed.
But this does not hold for any choice of quasi-distance, as 
demonstrated by the following example.

This example illustrates that the primal regularization 
may not help with the approximation of interpolation data $\Sigma$
that correspond to the theoretical data of 
some valid \psd outside the class of \psd \ used for interpolating.
The reason is that the quadratic term is increasing 
for large entries of $\Lambda$, but  not necessarily 
when approaching the boundary.

\begin{example}
{\sl
Consider now the approximative interpolation problem
$(\cP_\approx^2)$ for the special case that 
$d(\Phi||\Psi)=\frac{1}{2}\frac{(\Phi-\Psi)^2}{\Psi}$:


Form the Lagrangian function
\[
L(\Phi, \Delta;\q) \triangleeq \frac{1}{2} \int \frac{(\Phi-\Psi)^2}{\Psi}
+ \alpha \|\Delta\|^2  \qquad \quad \]
\[ \qquad
+ \sum_{k=0}^n  q_k\left( \Delta_k+  r_k \right)
-\ip{\Phi}{Q} 
\]

The optimal $\hat \Delta = -\frac{1}{2\alpha}q$.
The optimal \psd \ $\Phi$ is  determined 
by 
\[  \int \left( \frac{\Phi}{\Psi}-1-Q\right)\delta\Phi 
\, d\theta=0, \]

for all $\delta\Phi$, {\it i.e.} $\Phi = \Psi (Q+1)$.
The dual objective function is then given by 
\begin{equation}
\Omega (Q;\Psi)
= -\frac{1}{2}\ip{\Psi+\frac{1}{4\alpha}}{(Q+1)^2}+
 \sum_{k=0}^n  q_k r_k + \mbox{const.}
\end{equation}
Therefore, the regularization term and $\alpha$ only 
changes the prior and no matter how small $\alpha$ 
is chosen it is not always possible to find an 
interior point solution satisfying 
the stationarity conditions 
\[ r_k-\ip{\Psi(Q+1)}{z^k}=\ip{\frac{Q}{2\alpha}}{z^k}, \] 
for $k=0,1,\cdots,n$.
}\hfill$ \Box$
\end{example}

The next example illustrates that the dual regularization 
may not help with the approximation of interpolation data $\Sigma$
that does not correspond to the theoretical data of 
some valid \psd.
The reason is that the barrier function is increasing 
when approaching the boundary, but not necessarily 
for large entries of $\Lambda$.

\begin{example}
{\sl
Consider now the approximative interpolation problem
$(\cD_\approx^1)$ for the special case that
$\Omega_0 = -\trace\{\Lambda \Sigma \}+ \int \Psi \log Q$,
which corresponds to the primal with
the Kullback-Leibler divergence 
$d(\Phi||\Psi)=\Psi \log \frac{\Psi}{\Phi}$, 
and 
$B(Q)= \int \log Q$.

The objective function is then 
$ -\trace\{\Lambda \Sigma \}+ \int (\Psi+\lambda) \log Q$,
which corresponds to the exact interpolation problem 
with prior $\Psi+\lambda$.
If $\Sigma$ is not a positive semidefinite matrix, 
no matter how large $\lambda$ is, there exists 
no such exact interpolants, and the optimization 
problem $(\cD_\approx^1)$ has no finite optimum.
}\hfill$ \Box$
\end{example}

\bibliographystyle{IEEEbib}
\bibliography{bibliography_db}

\begin{thebibliography}{1}

\bibitem{GL}
T.T. Georgiou and A.~Lindquist,
\newblock ``Kullback-{L}eibler approximation of spectral density functions,''
\newblock {\em IEEE Transactions on Information Theory}, vol. 49, pp.
  2910--2917, Nov 2003.

\bibitem{HellingerMIMO}
A.~Ferrante, M.~Pavon, and F.~Ramponi,
\newblock ``Hellinger versus {K}ullback-{L}eibler multivariable spectrum
  approximation,''
\newblock {\em IEEE Trans. Automatic Control}, vol. 53, no. 4, pp. 954--967,
  May 2008.

\bibitem{Shore}
J.~Shore,
\newblock ``Minimum cross-entropy spectral analysis,''
\newblock {\em IEEE Trans. Acoustics, Speech and Signal Processing}, vol. 29,
  no. 2, pp. 230--237, Apr 1981.

\bibitem{ShoreJohnson}
J.~Shore and R.~Johnson,
\newblock ``Properties of cross-entropy minimization,''
\newblock {\em IEEE Trans. Information Theory}, vol. 27, no. 4, pp. 472--482,
  Jul 1981.

\bibitem{EK}
P.~Enqvist and J.Karlsson,
\newblock ``Minimal {I}takura-{S}aito distance and covariance interpolation,''
\newblock 2008, Conference on Decision and Control.

\bibitem{PerECC09}
P.~Enqvist,
\newblock ``Covariance interpolation and geometry of power spectral
  densities,''
\newblock in {\em Proceeding ECC 2009}, 2009.

\bibitem{Porat}
B.~Porat,
\newblock {\em Digital Processing of Random Signals, Theory \& Methods},
\newblock Prentice Hall, 1994.

\bibitem{BLunc}
C.I. Byrnes and A.~Lindquist,
\newblock {\em New Trends in Nonlinear Dynamics and Control}, chapter The
  uncertain generalized moment problem with complexity constraint, pp.
  267--278,
\newblock Springer Verlag, 2003.

\bibitem{EAtech}
P.~Enqvist and E.~Avvent,
\newblock ``Approximative linear and logarithmic interpolation of spectra,''
\newblock Tech. {R}ep. TRITA-MAT 09 OS 02, KTH Mathematics, 2009,
\newblock ISSN 1401-2294.

\end{thebibliography}

\end{document}